\documentclass[11pt]{article}
\usepackage{amssymb,amsmath,latexsym,verbatim,epsf,graphicx} 

\def\Z{\mathbb{Z}}
\def\R{\mathbb{R}}
\def\e{\mathrm{e}} 
\def\N{{\overline N}}

\begin{document}
\setlength{\parindent}{0pt}
\setlength{\parskip}{0.4cm}
\bibliographystyle{amsplain} 

\newtheorem{conjecture}{Conjecture}

\begin{center}

\Large{\bf Some Experimental Results on the Frobenius Problem}
\footnote{Appeared in \emph{Experimental Mathematics} {\bf 12}, no.~3 (2003), 263--269. 
          \emph{Mathematical Reviews Subject Numbers:} 05A15, 11P21; 11Y16.} 

\normalsize{\sc Matthias Beck, David Einstein, and Shelemyahu Zacks}

\end{center}

\footnotesize {\bf Abstract.} We study the \emph{Frobenius problem}: given relatively prime positive integers $a_{1} , \dots , a_{d}$, 
find the largest value of $t$ (the \emph{Frobenius number}) such that
$ \sum_{k=1}^d m_{k} a_{k} = t $ has no  solution in nonnegative integers $ m_{ 1 } , \dots , m_{ d } $. 
Based on empirical data, we conjecture that except for some special cases the Frobenius number can be bounded from above 
by $ \sqrt{a_{1}a_{2}a_{3}}^{5/4} - a_1 - a_2 - a_3 $. 
\normalsize 

\vspace{1cm}


\section{Introduction} 

Given positive integers $ a_{1} , \dots , a_{d} $ with $ \gcd (a_{1}, \dots, a_{d}) = 1 $, 
we call an integer $t$ \emph{representable}
if there exist nonnegative integers $ m_{ 1 } , \dots , m_{ d } $ such that
  \[  t = \sum_{j=1}^{d} m_{j} a_{j} \ . \]
In this paper, we discuss the \emph{linear Diophantine problem of Frobenius}: namely, find
the largest integer which is not representable. We call this largest integer the 
\emph{Frobenius number} $ g ( a_{1} , \dots , a_{d} ) $; 
its study was initiated in the 19th century. 
For $d=2$, it is well known (most probably at least since Sylvester \cite{sylvester}) that 
  \begin{equation}\label{gfor2} g(a_{1}, a_{2}) = a_{1} a_{2} - a_1 - a_2 \end{equation} 
For $d>2$, all attempts for explicit formulas have proved elusive. 
Two excellent survey papers on the Frobenius problem are \cite{alfonsin} and \cite{selmer}. 

Our goal is to establish bounds for $ g ( a_{1} , \dots , a_{d} ) $. 
The literature on such bounds is vast---see, for example, \cite{bdr,brauershockley,davison,erdosgrahamfrob,selmer,vitek}. 
We focus on the first non-trivial case $d=3$; any bound for this case yields a general bound, 
as one can easily see that $ g ( a_{1} , \dots , a_{d} ) \leq g ( a_{1} , a_2 , a_{3} ) $. 
All upper bounds in the literature are proportional to the product of two of the $a_k$'s. 
On the other hand, Davison proved in \cite{davison} the lower bound $ g ( a_{1} , a_2 , a_{3} ) \geq \sqrt{ 3 a_1 a_2 a_3 } - a_1 - a_2 - a_3 $. 
Experimental data (see Figure \ref{scat} below) shows that this bound is sharp in the sense 
that it is very often very close to $ g ( a_{1} , a_2 , a_{3} ) $. 
This motivates the question whether one can establish an \emph{upper} bound proportional 
to $ \sqrt{ a_1 a_2 a_3 }^{ p } $ where $p < 4/3 $. ($p = 4/3$ would be comparable to the known bounds.) 
In this paper we illustrate empirically, on the basis of more than ten thousand randomly chosen 
points, that $g(a_{1},a_{2},a_{3}) \leq \sqrt{a_{1}a_{2}a_{3}}^{5/4}-a_{1}-a_{2}-a_{3}.$ 


\section{Some geometric-combinatorial ingredients} 

Another motivation for the search for an upper bound proportional to $ \sqrt{ a_1 a_2 a_3 }^p $ comes from the following formula of \cite{bdr}, 
which is the basis for our study: 
Let $ a, b, c $ be pairwise relatively prime positive integers, and define 
  \[ N_{t} (a, b, c) := \# \left\{ (m_{1}, m_2, m_{3}) \in \Z^{3} : \ m_k \geq 0 , \ a m_1 + b m_2 + c m_3 = t \right\} \ . \] 
Then 
  \begin{eqnarray} 
    &\mbox{}& N_{t} (a, b, c) = \frac{ t^{ 2 }  }{ 2abc } + \frac{ t }{ 2 } \left( \frac{ 1 }{ ab } + \frac{ 1 }{ ac } + \frac{ 1 }{ bc }  \right) + \frac{ 1 }{ 12 } \left( \frac{ 3 }{ a } + \frac{ 3 }{ b } + \frac{ 3 }{ c } + \frac{ a }{ bc } + \frac{ b }{ ac } + \frac{ c }{ ab }  \right) \nonumber \\ 
    &\mbox{}& \qquad \qquad \qquad \qquad + \sigma_{-t} ( b, c; a ) + \sigma_{-t} ( a, c; b ) + \sigma_{-t} ( a, b; c ) \ , \label{N_t} 
  \end{eqnarray} 
where 
  \[ \sigma_t ( a, b; c ) := \frac{1}{c} \sum_{ \lambda^{ c } = 1 \not= \lambda } \frac{ \lambda^{ t } }{ \left( \lambda^{ a } - 1 \right) \left( \lambda^{ b } - 1 \right) }  \] 
is a \emph{Fourier-Dedekind sum}. 
One interpretation of $ N_{t} (a,b,c) $ is the number of partitions of $t$ with parts in the set $ \{ a,b,c \} $. 
Geometrically, $ N_{t} (a,b,c) $ enumerates integer points on the triangle 
  \[ \left\{ (x_{1}, x_2, x_{3}) \in \R^{3} : \ x_k \geq 0 , \ a x_1 + b x_2 + c x_3 = 1 \right\} \ , \] 
dilated by $t$. The Frobenius problem hence asks for the largest integer dilate of this triangle that contains no integer point, 
in other words, the largest $t$ for which $ N_{t} (a,b,c) = 0 $. 
It is also worth mentioning that the condition that $a$, $b$, and $c$ are \emph{pairwise} relatively prime is no restriction, due 
to Johnson's formula \cite{johnson}: if $ m = \gcd (a,b) $ then 
  \begin{equation}\label{johnsonid} g(a,b,c) = m \ g \left( \frac{a}{m} , \frac{b}{m} , c \right) + (m-1) c \ . \end{equation} 
In \cite{bdr}, formulas analogous to (\ref{N_t}) for $ d > 3 $ are given. In our case ($d=3$), a straightforward calculation shows 
  \[ \sigma_t ( a, b; c ) = - \frac{1}{4c} \sum_{k=1}^{c-1} \frac{ \e^{ \frac{\pi i k}{c} (-2t+a+b) }  }{ \sin \frac{ \pi k a }{ c } \sin \frac{ \pi k b }{ c } } \ . \] 
In fact, $ \sigma_t ( a, b; c ) $ is a \emph{Dedekind-Rademacher sum} \cite{rademacherdedekind}, as shown in \cite{bdr}. 
Hence we can rewrite (\ref{N_t}) as   
  \begin{eqnarray} 
    &\mbox{}& N_{t} (a,b,c) = \frac{ t^{ 2 }  }{ 2abc } + \frac{ t }{ 2 } \left( \frac{ 1 }{ ab } + \frac{ 1 }{ ac } + \frac{ 1 }{ bc }  \right) + \frac{ 1 }{ 12 } \left( \frac{ 3 }{ a } + \frac{ 3 }{ b } + \frac{ 3 }{ c } + \frac{ a }{ bc } + \frac{ b }{ ac } + \frac{ c }{ ab }  \right) \nonumber \\ 
    &\mbox{}& \qquad \qquad \qquad - \frac{1}{4a} \sum_{k=1}^{a-1} \frac{ \e^{ \frac{\pi i k}{a} (2t+b+c) }  }{ \sin \frac{ \pi k b }{ a } \sin \frac{ \pi k c }{ a } } - \frac{1}{4b} \sum_{k=1}^{b-1} \frac{ \e^{ \frac{\pi i k}{b} (2t+a+c) }  }{ \sin \frac{ \pi k a }{ b } \sin \frac{ \pi k c }{ b } } - \frac{1}{4c} \sum_{k=1}^{c-1} \frac{ \e^{ \frac{\pi i k}{c} (2t+a+b) }  }{ \sin \frac{ \pi k a }{ c } \sin \frac{ \pi k b }{ c } }   \ . \label{N_t2} 
  \end{eqnarray} 
If we write the ``periodic part'' of $ N_{t} (a,b,c) $ as 
  \[ P_{t} (a,b,c) := \frac{1}{4a} \sum_{k=1}^{a-1} \frac{ \e^{ \frac{\pi i k}{a} (2t+b+c) }  }{ \sin \frac{ \pi k b }{ a } \sin \frac{ \pi k c }{ a } } + \frac{1}{4b} \sum_{k=1}^{b-1} \frac{ \e^{ \frac{\pi i k}{b} (2t+a+c) }  }{ \sin \frac{ \pi k a }{ b } \sin \frac{ \pi k c }{ b } } + \frac{1}{4c} \sum_{k=1}^{c-1} \frac{ \e^{ \frac{\pi i k}{c} (2t+a+b) }  }{ \sin \frac{ \pi k a }{ c } \sin \frac{ \pi k b }{ c } }  \ , \] 
(\ref{N_t2}) becomes 
  \[ N_{t} (a,b,c) = \frac{ t^{ 2 }  }{ 2abc } + \frac{ t }{ 2 } \left( \frac{ 1 }{ ab } + \frac{ 1 }{ ac } + \frac{ 1 }{ bc }  \right) + \frac{ 1 }{ 12 } \left( \frac{ 3 }{ a } + \frac{ 3 }{ b } + \frac{ 3 }{ c } + \frac{ a }{ bc } + \frac{ b }{ ac } + \frac{ c }{ ab }  \right) - P_{t} (a,b,c) \ . \] 
If we can bound $ P_{t} (a,b,c) $ from above by, say, $B$ then the roots of $ N_{t} (a,b,c) $---and hence $ g (a,b,c) $---can be bounded from above: 
  \begin{eqnarray*} 
    &\mbox{}& g (a,b,c) \leq abc \left( - \tfrac{1}{2} \left( \tfrac{ 1 }{ ab } + \tfrac{ 1 }{ ac } + \tfrac{ 1 }{ bc } \right) + \sqrt{ \tfrac{1}{4} \left( \tfrac{ 1 }{ ab } + \tfrac{ 1 }{ ac } + \tfrac{ 1 }{ bc } \right)^2 - \tfrac{2}{abc} \left( \tfrac{ 1 }{ 12 } \left( \tfrac{ 3 }{ a } + \tfrac{ 3 }{ b } + \tfrac{ 3 }{ c } + \tfrac{ a }{ bc } + \tfrac{ b }{ ac } + \tfrac{ c }{ ab }  \right) - B \right) } \right) \\ 
    &\mbox{}& \qquad = - \tfrac{1}{2} ( a + b + c ) + \sqrt{ \tfrac{1}{4} (abc)^2 \left( \tfrac{ 1 }{ ab } + \tfrac{ 1 }{ ac } + \tfrac{ 1 }{ bc } \right)^2 - \tfrac{1}{6} abc \left( \tfrac{ 3 }{ a } + \tfrac{ 3 }{ b } + \tfrac{ 3 }{ c } + \tfrac{ a }{ bc } + \tfrac{ b }{ ac } + \tfrac{ c }{ ab }  \right) + 2 B \, abc } \\ 
    &\mbox{}& \qquad = \sqrt{ 2 B \, abc + \tfrac{1}{12} \left( a^2 + b^2 + c^2 \right) } - \tfrac{1}{2} ( a + b + c ) 
  \end{eqnarray*} 
From this computation, the question of the existence of an upper bound for $ g ( a,b,c ) $ proportional to $ \sqrt{ abc }^p $ comes up naturally. 
Unfortunately, it is not clear how to bound the periodic part $ P_{t} (a,b,c) $ effectively. An almost trivial bound for $ P_{t} (a,b,c) $ 
yielded in \cite{bdr} the inequality 
  \[ g(a,b,c) \leq \frac{ 1 }{ 2 } \left( \sqrt{ abc \left( a + b + c \right) } - a - b - c \right) \ , \] 
which is of comparable size to the other upper bounds for $ g ( a,b,c ) $ in the literature. 
However, we believe one can obtain bounds of smaller magnitude. 


\section{Special cases}

On the path to such ``better'' bounds, we first have to exclude some cases which definitely yield Frobenius numbers of size $a_k^2$. 
One of these cases are triples $(a,b,c)$ such that $c$ is representable by $a$ and $b$: by (\ref{gfor2}), we obtain in this 
case $ g(a,b,c) = ab - a - b $.  

A second case of triples $(a,b,c)$ that we need to exclude are those for which $a|(b+c)$. Brauer and Shockley \cite{brauershockley} 
proved that in this case 
  \[ g(a,b,c) = \max \left( b \left\lfloor \frac{ ac }{ b+c } \right\rfloor , c \left\lfloor \frac{ ab }{ b+c } \right\rfloor \right) - a \ . \] 
Here $ \left\lfloor x \right\rfloor $ denotes the greatest integer not exceeding $x$. 

An even less trivial example of special cases was given by Lewin \cite{lewin}, who studied the Frobenius number of \emph{almost arithmetic sequences}: 
If $m,n>0$, gcd$(a,n)=1$, and $ d \leq a $, then 
  \[ g \left( a, ma+n, ma+2n, \dots , ma + \left( d-1 \right) n \right) = \left( m \left\lfloor \frac{ a-2 }{ d-1 } \right\rfloor + m - 1 \right) a + \left( a-1 \right) n \ . \] 
For arithmetic sequences ($m=1$), this formula goes back to Roberts \cite{roberts}, for consecutive numbers ($m=n=1$) it is due to Brauer \cite{brauer}. 
For the special case $d=3$, we obtain 
  \[ g ( a, ma+n, ma+2n ) = \left( m \left\lfloor \frac{a}{2} \right\rfloor - 1 \right) a + \left( a-1 \right) n \ .  \] 
As a function in $a, \ b:= ma+n, \ c:= ma + 2n$, this Frobenius number grows proportionally to $ab$, which means an upper bound proportional to 
$ \sqrt{abc}^p $ with $ p < 4/3 $ can not be achieved. 
Hence in our computations and conjectures about upper bounds for $ g(a,b,c) $, we will exclude the cases of one of the numbers being representable 
by the other two, one number dividing the sum of the other two, 
and almost arithmetic sequences. Finally, as noted above, thanks to (\ref{johnsonid}) we may assume without loss of generality 
that $a$, $b$, and $c$ are pairwise coprime. The triples $(a,b,c)$ that are not excluded will be called \emph{admissible}. 


\section{Computations} 

In the present section we discuss the computation of the Frobenius number. 
For convenience we computed the number $ f(a,b,c) = g(a,b,c) + a + b + c $. 
It is not hard to see that $f(a,b,c)$ is the largest integer that can not be 
represented by a linear combination of $a$, $b$, and $c$ with \emph{positive} integer 
coefficients. The respective counting function 
  \[ \N_{t} (a, b, c) := \# \left\{ (m_{1}, m_2, m_{3}) \in \Z^{3} : \ m_k > 0 , \ a m_1 + b m_2 + c m_3 = t \right\} \ . \] 
can also be found in \cite{bdr} and is closely related to $ N_{t} (a,b,c) $: 
  \begin{eqnarray*} 
    &\mbox{}& \N_{t} (a,b,c) = \frac{ t^{ 2 }  }{ 2abc } - \frac{ t }{ 2 } \left( \frac{ 1 }{ ab } + \frac{ 1 }{ ac } + \frac{ 1 }{ bc }  \right) + \frac{ 1 }{ 12 } \left( \frac{ 3 }{ a } + \frac{ 3 }{ b } + \frac{ 3 }{ c } + \frac{ a }{ bc } + \frac{ b }{ ac } + \frac{ c }{ ab }  \right) \\ 
    &\mbox{}& \qquad \qquad \qquad \qquad - \frac{1}{4a} \sum_{k=1}^{a-1} \frac{ \e^{ \frac{\pi i k}{a} (-2t+b+c) }  }{ \sin \frac{ \pi k b }{ a } \sin \frac{ \pi k c }{ a } } - \frac{1}{4b} \sum_{k=1}^{b-1} \frac{ \e^{ \frac{\pi i k}{b} (-2t+a+c) }  }{ \sin \frac{ \pi k a }{ b } \sin \frac{ \pi k c }{ b } } - \frac{1}{4c} \sum_{k=1}^{c-1} \frac{ \e^{ \frac{\pi i k}{c} (-2t+a+b) }  }{ \sin \frac{ \pi k a }{ c } \sin \frac{ \pi k b }{ c } }   \ . 
  \end{eqnarray*} 
The following illustrates our algorithm. 

\begin{verbatim} 
STEP 0: Initiate the intervals I1, I2, I3 for the selection of the 
        arguments a,b,c; 

STEP 1: Draw at random integers a,b,c from I1, I2, I3 respectively;

STEP 2: Test a,b,c for coprimality and for almost arithmetic sequences;

STEP 3: IF (a,b,c are not pairwise coprime) or IF (a,b,c are almost arithmetic)
        {discard a,b,c and GOTO STEP 1}
        ELSE {SET delta <- min(a,b,c); GOTO STEP 4};

STEP 4: Compute z=sqrt(a*b*c), 
        SET mb <- INT(sqrt(3)*z)+delta; 
        SET t <- mb;

STEP 5: Compute NB(t,a,b,c);

STEP 6: IF (NB(t,a,b,c)>0) {SET t <- t-1, and GOTO STEP 5}
        ELSE {GOTO STEP 7};

STEP 7: SET f <- t;

STEP 8: IF(mb-f < delta) {SET 
        mb <- mb+delta
        t <- mb
        GOTO STEP 5}
        ELSE {GOTO STEP 9};

STEP 9: PRINT f(a,b,c) <- f;
        STOP.
\end{verbatim} 

For example, for $a=7, b=13 , c=30$ the program yields the Frobenius
number $ f(7,13,30) = 95 $, or $ g(7,13,30) = 45 $. This program was tested against arguments which yield known
results, and found to be correct.

Our program is to choose at random arguments $a, b, c$ in a certain range (in our case $[1,750]$), and test the 
triplets for admissibility. For admissible triplets $a,b,c$ we compute the Frobenius number $f(a,b,c)$ based on 
the straightforward observation that, once we have $a=\min(a,b,c)$ 
consecutive integers which are representable, we know that every integer beyond that interval is representable as well. 
We start searching for roots of $\N_{t} (a,b,c)$ at the lower bound $\sqrt{3abc}$. If a root is found at an integer 
$f$, we repeat this search until we found an interval of $a$ integers $t$ with $\N_{t} (a,b,c) > 0$, that is, an 
interval of $a$ representable integers. At this stopping point, the integer $f$ is the sought-after Frobenius number $f(a,b,c)$. 

We have created a {\tt PARI-GP} program\footnote{Our program can be downloaded at {\tt www.math.binghamton.edu/matthias/frobcomp.html}}, 
following the above algorithm. The program proved to be quite efficient, since most of the 
values of $f(a,b,c)$ were found to be close to the lower bound $\sqrt{3abc}$, 
as shown in the analysis below. The Dedekind-Rademacher sums appearing in (\ref{N_t2}) can be computed 
very efficiently because they satisfy a reciprocity law 
(\cite{rademacherdedekind}, for computational complexity see also \cite{knuth}), 
which allows us to calculate their values similar in spirit 
to the Euclidean algorithm. This implies that for a given $t$, $\N_{t} (a,b,c)$ can be computed with our rather simple algorithm in $O(\log(c))$ time, assuming that 
$c=\max(a,b,c)$. Hence if $f(a,b,c)$ is close to the lower bound $\sqrt{3abc}$---which, again, happens in the vast majority 
of cases---, we obtain $f(a,b,c)$ in $O(a\log(c))$ time. On the other hand, we can of course not assume that $f(a,b,c)$ is close to $\sqrt{3abc}$; 
still we get `at worst' a computation time of $O(ab\log(c))$. 
What makes this analysis even more appealing is that it applies to the general case of the Frobenius problem. 
As mentioned above, there is an analog for (\ref{N_t}) and (\ref{N_t2}) for $ d > 3 $ \cite{bdr}, which again 
is a lattice-point count in a polytope and as such is known (for fixed $d$) to be computable in $ O \left( p \left( \log a_1, \dots, \log a_d \right) \right) $ 
time for some polynomial $p$ \cite{barvinokehrhart}. With an analogous algorithm for the general case, we would hence be able to compute $f(a_1, \dots, a_d)$ in 
$ O \left( a_1 a_2 \, p \left( \log a_1, \dots, \log a_d \right) \right) $ time, where $a_1 < a_2 < \dots < a_d$. 
As in the three-variable case---in fact, even more so---, most Frobenius numbers will be situated very close to the lower bound $\sqrt{3 a_1 a_2 a_3}$, which 
means that in the vast majority of cases we can expect a computation time of $ O \left( a_1 \, p \left( \log a_1, \dots, \log a_d \right) \right) $. 

The computational complexity of the Frobenius problem is very interesting and still gives rise to ongoing studies. 
Davison \cite{davison} provided an algorithm for the three-variable case ($a<b<c$) which runs in $O(\log b)$ time. 
The general case is still open. While Kannan \cite{kannan} proved that there is a polynomial-time algorithm (polynomial 
in $\log a_1, \dots, \log a_d$) to find $ g ( a_{1} , \dots , a_{d} ) $ for fixed $d$, no such algorithm is known for 
$d>3$. The fastest general algorithm which we are aware of is due to Nijenhuis \cite{nijenhuis} and runs in 
$O( d \, a \log a ) $ time, where $a = \min ( a_{1} , \dots , a_{d} ) $. 
Hence, while our primitive algorithm is not competitive for the three-variable case of the Frobenius problem, it might be worthwhile 
to develop it further in the general case. 

We initially implemented our program as an {\tt MS-DOS QUICK BASIC} program and experienced 
some interesting problems due to floating-point errors: computing generalized Dedekind sums 
can get challenging for large arguments. These problems were only discovered when we reimplemented 
the algorithm in {\tt PARI-GP}, which has an extended precision aritmetic and also keeps track of
roundoff errors effectively. It is worth mentioning that both Knuth's algorithm \cite{knuth} for 
the computation of Rademacher-Dedekind sums and Davison's algorithm \cite{davison} for computing 
$g(a,b,c)$ are integer algorithms and therefore are very stable.

With our program we generated at random 10000 admissible 
triplets. Our main question is the relation of the Frobenius number $f = f(a,b,c)$ to $z:=\sqrt{abc}$. 
The following is a statistical description of the ratios $R:=f/z$. 

\textbf{Descriptive Statistics} \\ 
\begin{tabular}{ccccccccc} Variable & $N$ & Mean & Median & StDev & Min & Max & $Q_1$ & $Q_3$ \\ 
                           $R$ & 10000 & 2.283 & 2.012 & 0.737 & 1.736 & 9.332 & 1.940 & 2.299 \end{tabular} 

$ Q_1 $ and $ Q_3 $ are the first and third quartiles, respectively. 
We see in the above table that 50\% of the cases have a ratio smaller than
2.01, and 75\% have ratio smaller than 2.30. In the following figures we
present a box-plot and a histogram of the variable $R$. 

\begin{figure}[hb!]
\begin{center}
\includegraphics[totalheight=3in]{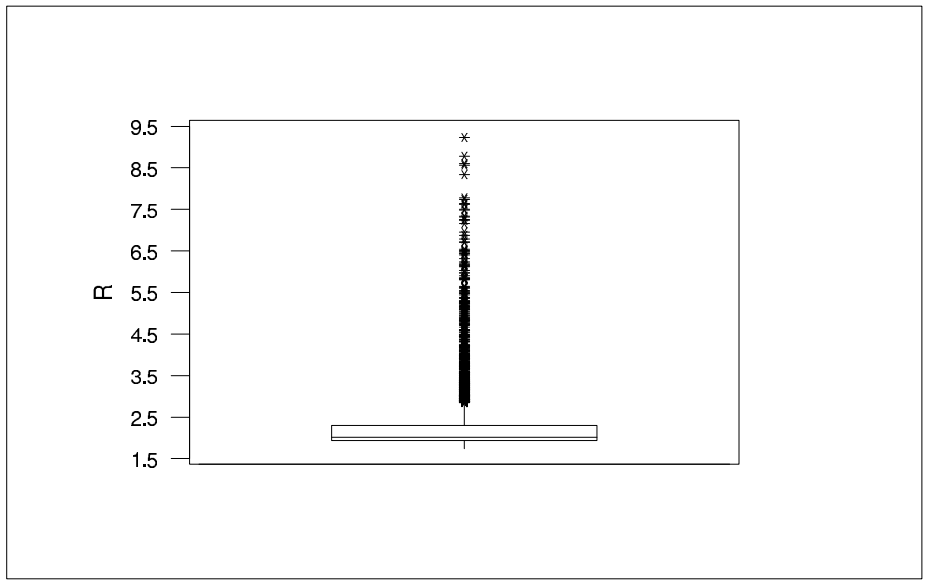} 
\end{center}
\caption{Box-plot} \label{boxplot} 
\end{figure} 

\begin{figure}[hbt]
\begin{center}
\includegraphics[angle=270,totalheight=3in]{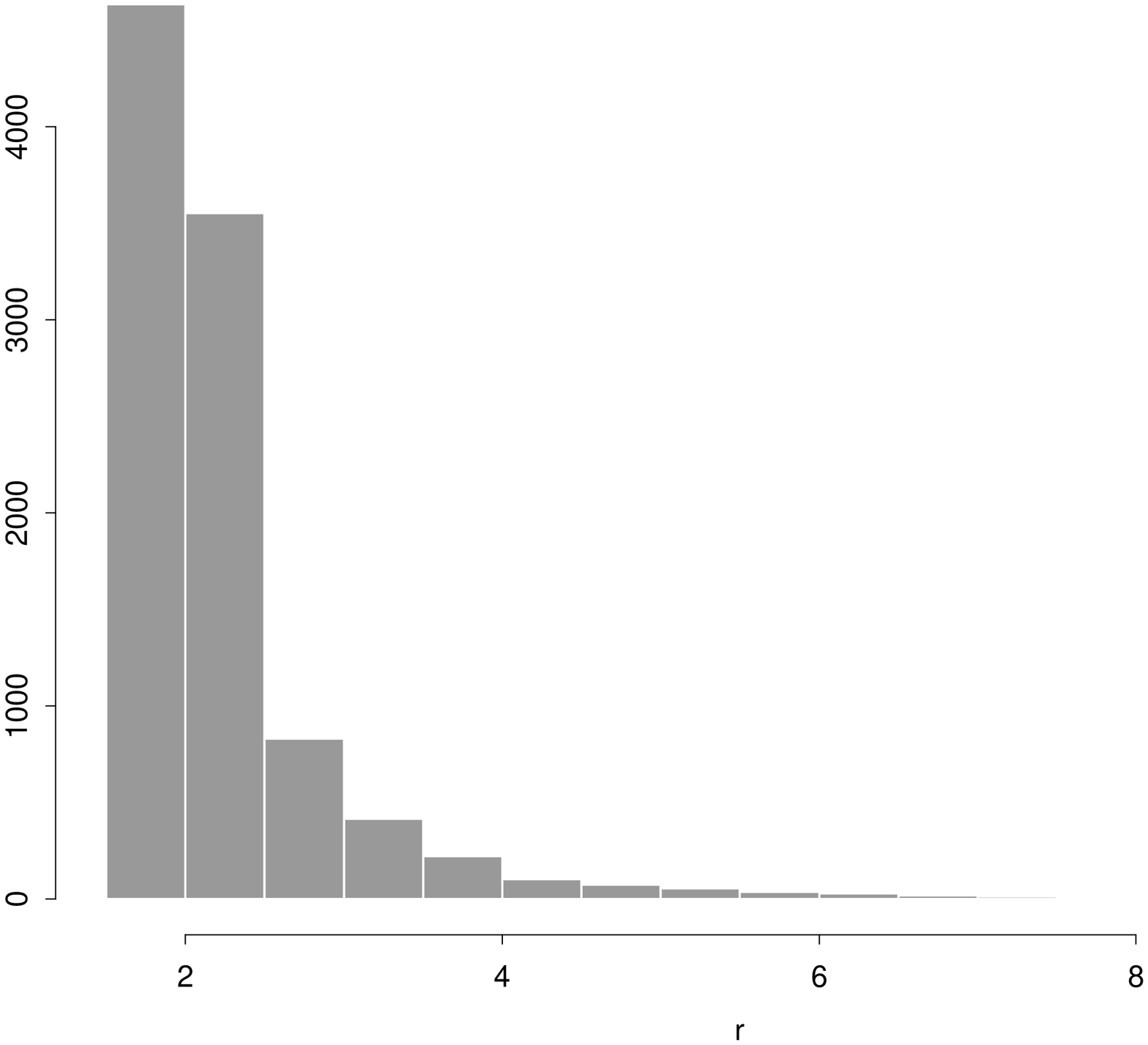} 
\end{center}
\caption{Histogram} 
\end{figure} 

In the box-plot, the bottom line of the box corresponds to the first
quartile $Q_{1}$. The top line of the box corresponds to the third quartile $Q_{3}$. 
There are 980 points above the value of $R=3$. Only 24 points, which are listed below, have a
value of $R$ greater than 7.5. 

\begin{center}
\begin{tabular}{r|r|r|r|r|r|r|r} 
$a \ $ & $b \ $ & $c \ $ & $f(a,b,c)$ & $z=\sqrt{abc}$ & $\sqrt 3 \, z \ $ & $z^{5/4} \ $ & $R=f/z$ \\ \hline 
487 & 733 & 738 & 121755 & 16231.0 & 28112.9 & 183202 & 7.50140 \\ 
229 & 483 & 662 & 64901 & 8557.0 & 14821.1 & 82300 & 7.58457 \\ 
223 & 307 & 698 & 52657 & 6912.7 & 11973.2 & 63032 & 7.61740 \\ 
244 & 357 & 619 & 56067 & 7343.0 & 12718.5 & 67974 & 7.63542 \\ 
509 & 541 & 557 & 95788 & 12384.7 & 21450.9 & 130649 & 7.73439 \\ 
262 & 349 & 699 & 61861 & 7994.7 & 13847.2 & 75597 & 7.73776 \\ 
475 & 611 & 679 & 109183 & 14037.9 & 24314.4 & 152802 & 7.77773 \\ 
248 & 305 & 439 & 45274 & 5762.5 & 9980.9 & 50207 & 7.85671 \\ 
265 & 488 & 509 & 65434 & 8113.2 & 14052.5 & 77000 & 8.06514 \\ 
274 & 401 & 695 & 70596 & 8738.6 & 15135.6 & 84489 & 8.07868 \\ 
368 & 415 & 599 & 77374 & 9564.5 & 16566.2 & 94586 & 8.08972 \\ 
281 & 341 & 502 & 57790 & 6935.6 & 12012.8 & 63293 & 8.33241 \\ 
315 & 488 & 559 & 77734 & 9269.8 & 16055.8 & 90958 & 8.38571 \\ 
305 & 319 & 652 & 67142 & 7964.7 & 13795.3 & 75242 & 8.42995 \\ 
393 & 452 & 619 & 89830 & 10486.0 & 18162.3 & 106112 & 8.56664 \\ 
313 & 532 & 579 & 84150 & 9819.0 & 17007.0 & 97743 & 8.57012 \\ 
301 & 479 & 725 & 87903 & 10224.0 & 17708.5 & 102808 & 8.59773 \\ 
655 & 671 & 679 & 150043 & 17274.9 & 29921.1 & 198048 & 8.68558 \\ 
296 & 731 & 749 & 110834 & 12730.5 & 22049.9 & 135225 & 8.70618 \\ 
359 & 520 & 619 & 94318 & 10749.6 & 18618.9 & 109457 & 8.77406 \\ 
337 & 346 & 701 & 79559 & 9040.9 & 15659.3 & 88159 & 8.79989 \\ 
320 & 469 & 491 & 77556 & 8584.2 & 14868.4 & 82628 & 9.03469 \\ 
335 & 668 & 669 & 112894 & 12235.6 & 21192.6 & 128685 & 9.22672 \\ 
379 & 389 & 748 & 97998 & 10501.4 & 18188.9 & 106306 & 9.33194 
\end{tabular} 
\end{center} 

\begin{figure}[htb]
\begin{center}
\includegraphics[angle=-90,totalheight=5in]{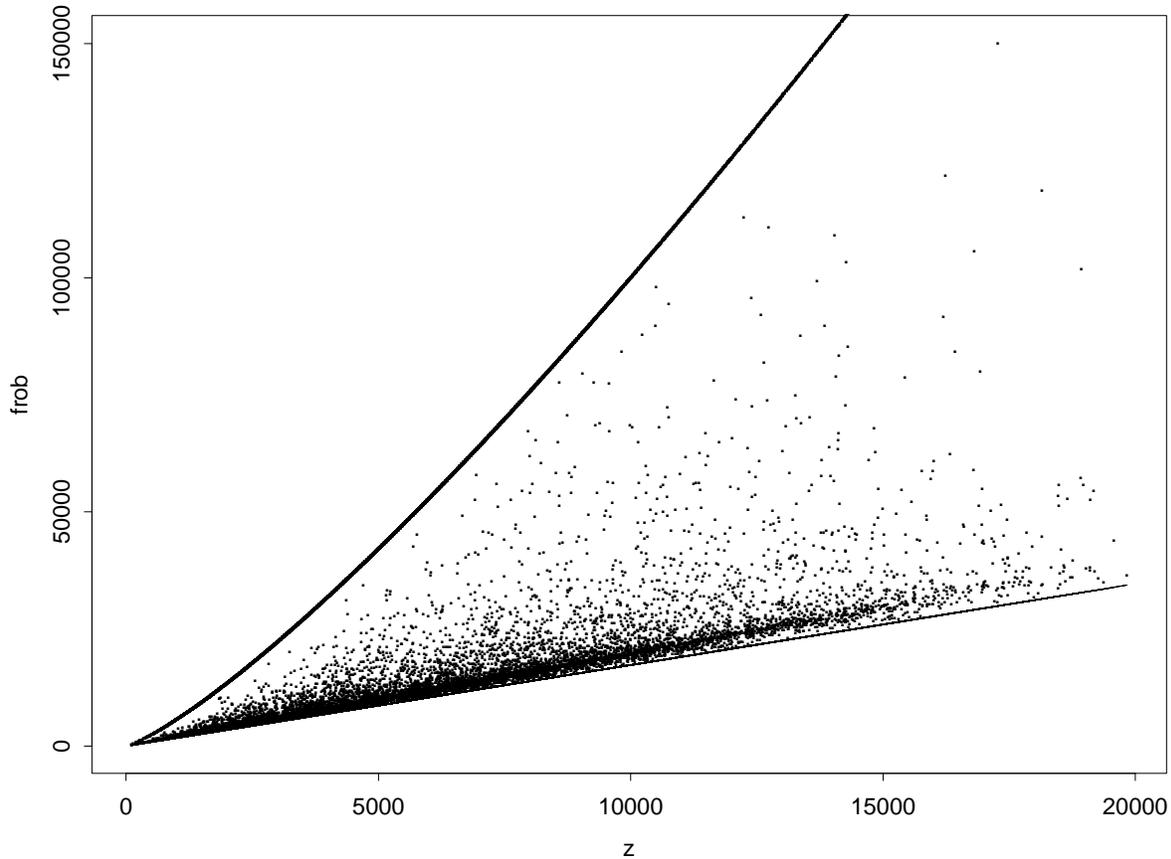}
\end{center}
\caption{$f=f(a,b,c)$ as a function of $ z = \sqrt{abc} $}\label{scat} 
\end{figure}

In Figure \ref{scat} we presents all the points $(z,f)$. 
Notice that all the points are above the straight line $\sqrt{3}z$, which illustrates Davison's lower bound. 
The upper bound is, however, convex. It is included in the figure as the graph $z^{5/4}$. 



\section{Conjectures and closing remarks} 

Randomly chosen admissible arguments tend to yield a Frobenius number $f$ smaller than the expected number (mean) which is 
estimated to be $ 2.28 z $. The distribution of $R=f/z$ is very skewed (positive asymmetry) as seen in the Figure \ref{boxplot}. Since 
10000 random points yielded $ f < z^{5/4} $, or $ g(a,b,c) < \sqrt{abc}^{5/4}-a-b-c $, the probability that a future randomly 
chosen admissible triplet with $ z = \sqrt{abc} < 20000 $ will yield $ f \geq z^{5/4} $ is smaller than $ 1/10000 $. 

In general, our data suggests that one can obtain an upper bound of smaller magnitude than what the above cited results state. 
Again, the upper bounds in the literature are comparable to an upper bound proportional to $ \sqrt{ abc }^{ 4/3 } $. 
We believe the following is true. 
\begin{conjecture} There exists an upper bound for $ (a, b, c) $ proportional to $ \sqrt{ abc }^{ p } $ where $ p < \frac{4}{3} $, valid for all 
admissible triplets $(a,b,c)$. \end{conjecture} 
In fact, our data suggests, more precisely, that for all admissible triplets $(a,b,c)$, 
  \[ g(a, b, c) \leq \sqrt{abc}^{5/4}-a-b-c \ . \] 
It is very improbable that a randomly chosen admissible triple $(a,b,c)$, such that $ \sqrt{abc} < 20000 $, 
will yield $g(a, b, c) > \sqrt{abc}^{5/4}-a-b-c$. However, we remark that there might be specific structures 
of triples $(a,b,c)$, close to almost arithmetic, for which $g(a, b, c) > \sqrt{abc}^{5/4}-a-b-c$. 
This is generally not the case. 

\vspace{1cm} 
{\bf Acknowledgements}. We would like to thank Tendai Chitewere for days and days of computing time, 
Gary Greenfield for the non-trivial task of converting our pictures into \LaTeX-friendly postscript, 
and the referee and associate editor for helpful comments on the first version of this paper. 
Finally, we would like to thank the authors and maintainers of {\tt PARI-GP}.


\bibliographystyle{plain}

\def\cprime{$'$}
\providecommand{\bysame}{\leavevmode\hbox to3em{\hrulefill}\thinspace}
\providecommand{\MR}{\relax\ifhmode\unskip\space\fi MR }
\providecommand{\MRhref}[2]{%
  \href{http://www.ams.org/mathscinet-getitem?mr=#1}{#2}
}
\providecommand{\href}[2]{#2}

\sc Department of Mathematical Sciences\\
    State University of New York\\
    Binghamton, NY 13902-6000 \\ 
\tt matthias@math.binghamton.edu 

\sc Structured Decisions Corp.\\ 
    1105 Washington St.\\ 
    West Newton MA 02465\\ 
\tt deinst@world.std.com 

\sc Department of Mathematical Sciences\\
    State University of New York\\
    Binghamton, NY 13902-6000 \\ 
\tt shelly@math.binghamton.edu
    
\end{document}